\documentclass[letterpaper, 10 pt, conference]{ieeeconf}  % Comment this line out
\IEEEoverridecommandlockouts                              % This command is only
\usepackage[hyphenbreaks]{breakurl}
\usepackage{cite}
\usepackage{amsmath,amssymb}
\usepackage{rotating}
\usepackage{tablefootnote}
\usepackage{url}
\usepackage{hyperref}

\title{\LARGE \bf
Verifying Global Optimality of Candidate Solutions to Polynomial\\ Optimization Problems using a Determinant Relaxation Hierarchy
}

\author{Sikun Xu, Ruoyi Ma, Daniel K. Molzahn, Hassan Hijazi, and C\'edric Josz % <-this % stops a space
\thanks{This work was supported by NSF PD 18-7607 Energy, Power, Control, and Networks, award number 2023032 and 2023140.}% <-this % stops a space
\thanks{S. Xu, R. Ma, and C. Josz are with IEOR Columbia University, S. W. Mudd Building, 500 W 120th St, New York, NY 10027-6623 {\tt\small sx2265@columbia.edu, rm3715@columbia.edu, cj2638@columbia.edu}}%
\thanks{D. K. Molzahn is with the
School of Electrical and Computer Engineering, Georgia Institute of Technology, 
E176 Van Leer Building
777 Atlantic Dr. NW
Atlanta, GA 30313, USA
        {\tt\small molzahn@gatech.edu}}%
\thanks{H. Hijazi is with the Los Alamos National Laboratory
{\tt\small hlh@lanl.gov}}
}

\begin{document}

\maketitle
\pagestyle{plain}

%%%%%%%%%%%%%%%%%%%%%%%%%%%%%%%%%%%%%%%%%%%%%%%%%%%%%%%%%%%%%%%%%%%%%%%%%%%%%%%%
\begin{abstract}
We propose an approach for verifying that a given feasible point for a polynomial optimization problem is globally optimal. The approach relies on the Lasserre hierarchy and the result of Lasserre regarding the importance of the convexity of the feasible set as opposed to that of the individual constraints. By focusing solely on certifying global optimality and relaxing the Lasserre hierarchy using necessary conditions for positive semidefiniteness based on matrix determinants, the proposed method is implementable as a computationally tractable linear program. We demonstrate this method via application to several instances of polynomial optimization, including the optimal power flow problem used to operate electric power systems.
\end{abstract}

%%%%%%%%%%%%%%%%%%%%%%%%%%%%%%%%%%%%%%%%%%%%%%%%%%%%%%%%%%%%%%%%%%%%%%%%%%%%%%%%
\section{INTRODUCTION}

Consider a polynomial optimization problem
\begin{equation}
\label{eq:pop}
\begin{array}{ll}
\inf_x & f(x) := \sum_\alpha f_{\alpha} x^\alpha \\
\mathrm{s.t.} & g_i(x) := \sum_\alpha g_{i,\alpha} x^\alpha \geqslant 0, \quad i=1,\ldots,m
\end{array}
\end{equation}
where we use the multi-index notation $x^\alpha := x_1^{\alpha_1} \cdots x_n^{\alpha_n}$ for $x \in {\mathbb R}^n$,
$\alpha \in {\mathbb N}^n$, and
where the data are polynomials $f, g_1, \ldots, g_m \in {\mathbb R}[x]$
so that in the above sums only a finite number of coefficients
$f_{\alpha}$ and $g_{i,\alpha}$ are nonzero. We will use the notation $|\alpha|:= \sum_{k=1}^n \alpha_k$.

Polynomial optimization problems~\eqref{eq:pop} are generally non-convex and may therefore have multiple local and global optima. Local optima are feasible points for which no nearby point has a lower-cost objective value. Global optima are feasible points which have the lowest objective value among all feasible points. Many applications require (or at least benefit from) finding a global optimum to~\eqref{eq:pop}. Provably obtaining a global optimum is often accomplished by:
\begin{enumerate}
\item Using a local solver to compute a locally optimal point that is a \emph{candidate} for being a global optimum.
\item Computing the \emph{optimality gap} based on the difference between the objective value of the candidate point and the lower bound on the objective value provided by a convex relaxation of~\eqref{eq:pop}.
\end{enumerate}
The candidate point is declared to be globally optimal if the optimality gap is sufficiently small. When the optimality gap is not sufficiently small, typical global optimization methods attempt to close the optimality gap by applying a variety of techniques for tightening the relaxation and obtaining a better local optimum~\cite{floudas2013deterministic,GOPINATH2020106688}.

For certain polynomial optimization problems, especially those modeling physical systems, it is often straightforward to obtain high-quality local optima due to preexisting knowledge of appropriate initializations for local solvers.
%Despite including particularly challenging instances, knowledge of appropriate initializations often results in local solvers yielding globally optimal solutions to certain classes of practically relevant non-convex polynomial optimization problems. 
Consider, for instance, the optimal power flow (OPF) problem used in the operation of electric power systems~\cite{carpentier-1962,cain-2012}. OPF problems are known to be NP-hard~\cite{bienstock2015nphard,pascalNPhard}. Accordingly, there exist various test cases that challenge many solution algorithms~\cite{bukhsh2013,pglib-opf,kocuk2015,narimani_molzahn_wu_crow-empirical_nonconvexity_study}. However, despite this challenging worst-case complexity, local solvers with physically justified initializations find the global optima for many practical OPF problems as verified by the small or zero optimality gaps obtained via a variety of convex relaxation techniques~\cite{lavaei_tps,low_tutorial,coffrin2015,mh_sparse_msdp,josz2018lasserre,kocuk2017minor}. (See~\cite{molzahn_hiskens-fnt2017} for a detailed survey.) Moreover, appropriately initialized local solvers can be significantly more computationally tractable than the convex relaxations used to certify global optimality~\cite{ferc5,matpower}. Thus, the computational effort for provably obtaining a globally optimal solution is often dominated by the time spent solving relaxations which show that a candidate feasible point from a local solver is, in fact, globally optimal.

Such problems motivate the development of methods for checking whether a given feasible point (e.g., the output of a local solver) is globally optimal. In contrast to convex relaxation techniques, whose objective values provide bounds on the potential suboptimality associated with a local solution, this paper solely focuses on methods for \emph{certifying} whether a given feasible point is globally optimal. This less ambitious aim in comparison to existing convex relaxation techniques facilitates computational advantages for the proposed approach. Specifically, we reduce the problem of certifying global optimality of a candidate point to testing feasibility of a linear program.

Our proposed approach is particularly useful in settings that require the solution of many polynomial optimization problems. When applied as  a screening step, our approach limits the need for explicitly solving computationally intensive relaxations to problems for which the candidate point and the selected relaxation result in a non-zero optimality gap. This allows us to leverage much of the extensive theoretical work regarding relaxations of polynomial optimization problems that has heretofore been practically limited by computational challenges associated with existing conic optimization solvers.

We close the introduction by reviewing two papers that aim to certify global optimality of polynomial optimization problems. Reviews of the specific relaxation techniques we employ will be presented in the following section.

The approach in~\cite{lasserre2013inverse} finds a objective function $\hat{f}$ for which a candidate point, $\hat{x}$, is globally optimal:
\begin{equation}
\label{eq:closest_obj}
\begin{array}{ll}
\inf_{\hat{f}} &  \| f  - \hat{f} \| \\
\mathrm{s.t.} & \hat{f}(x) - \hat{f}(\hat{x}) \geqslant 0, ~ \forall x \in \mathcal{K}
\end{array}
\end{equation}
where $\mathcal{K}$ denotes the feasible space defined by the constraints in~\eqref{eq:pop} and \mbox{$\| \cdot \|$} is, for instance, the 2-norm of the coefficients. If the optimal objective value of~\eqref{eq:closest_obj} is equal to zero, then $\hat{x}$ is globally optimal. If not, then the solution of~\eqref{eq:closest_obj} provides insight into how far $\hat{x}$ is from the global optimum. Specifically, the solution to~\eqref{eq:closest_obj} yields an objective function, as close to the original objective function as possible, for which the candidate point $\hat{x}$ is globally optimal. In order to solve~\eqref{eq:closest_obj}, \cite{lasserre2013inverse} proposes to replace the non-negativity constraint using Putinar's Positivstellensatz. The drawback of this approach is that it is \textit{a priori} as computationally costly as actually computing a global solution.

In~\cite{molzahn_lesieutre_demarco-global_optimality_condition}, another related approach provides a sufficient (but not necessary) condition for certifying global optimality of a candidate point by evaluating whether the point satisfies the Karush-Kuhn-Tucker (KKT) conditions~\cite{karush1939,kt1951} for a relaxation. Assuming certain constraint qualification conditions (e.g., Slater's condition), a local solution necessarily satisfies the KKT conditions for the non-convex polynomial optimization problem~\eqref{eq:pop}. If the local solution additionally satisfies the KKT conditions for a convex relaxation of~\eqref{eq:pop}, then the local solution is, in fact, a global optimum. The approach in~\cite{molzahn_lesieutre_demarco-global_optimality_condition} specifically tests whether a local solution to an OPF problem satisfies the KKT conditions for the Shor relaxation~\cite{shor-1987b}, which is equivalent, in the case of OPF problems, to the first-order moment/sum-of-squares relaxation in the Lasserre hierarchy to be discussed below. While explicitly computing the Shor relaxation requires the solution of a semidefinite program, evaluating whether the primal/dual variables obtained from a local solver satisfy the Shor relaxation's KKT conditions only requires a single Cholesky factorization of a sparse matrix and a vector dot product. Thus, certifying global optimality by evaluating the condition in~\cite{molzahn_lesieutre_demarco-global_optimality_condition} is up to two orders of magnitude faster than explicitly solving the Shor relaxation for various OPF test cases, even after applying advanced techniques for improving the relaxation's computational tractability~\cite{jabr11,molzahn_holzer_lesieutre_demarco-large_scale_sdp_opf}. The approach proposed in the remainder of this paper can be viewed as an extension of the condition in~\cite{molzahn_lesieutre_demarco-global_optimality_condition} to exploit certain tighter relaxations, specifically, higher-order moment/sum-of-squares relaxations in the Lasserre hierarchy. 

The remainder of paper is organized as follows. Section~\ref{sec:approach} presents the proposed approach. Section~\ref{sec:Numerical experiments} applies the approach to several polynomial optimization problems including an instance of the optimal power flow problem.
%Section~\ref{sec:computational_details} describes details regarding the exploitation of sparsity and selective application of the higher-order relaxation constraints in order to obtain a tractable implementation.
%Section~\ref{sec:acopf} introduces the OPF problem and presents a small illustrative test case. 
Section~\ref{sec:conclusion} concludes the paper and proposes avenues for future research.

\section{Approach}
\label{sec:approach}

In 2001, the Lasserre hierarchy \cite{lasserre-2000,lasserre-2001} (see also \cite{parrilo-2000b,parrilo-2003}) was proposed to find global solutions to polynomial optimization problems. It is also known as moment/sum-of-squares hierarchy in reference to the primal moment hierarchy and the dual sum-of-squares hierarchy. Its global convergence is guaranteed by Putinar's Positivstellensatz \cite{putinar-1993} proven in 1993. Typically, if one of the constraints is a ball $x_1^2 + \hdots + x_n^2 \leqslant 1$, then the sequence of lower bounds provided by the hierarchy converges to the global infimum of the polynomial optimization problem. In addition, there is zero duality at all relaxation orders~\cite{josz-2015}.

The moment problem of order $d$ is defined as
\begin{equation}
\label{eq:lasserre}
\begin{array}{ll}
\inf_y & L_y(f)  \\
\mathrm{s.t.} & y_0 = 1 \\
 & M_d(y) \succcurlyeq 0 \\
 &  M_{d-k_i}(g_iy) \succcurlyeq 0, \quad i=1,\ldots,m
\end{array}
\end{equation}
where $\succcurlyeq 0$ denotes positive semidefiniteness, and where the Riesz functional, the moment matrix, and the localizing matrices are respectively defined by
\begin{equation}
\begin{array}{l}
L_y(f) := \sum_{\alpha} f_{\alpha} y_\alpha\\[.2cm]
M_d(y) := (y_{\alpha+\beta})_{|\alpha|,|\beta|\leqslant d} \\[.2cm]
M_{d-k_i}(g_iy) := (\sum_{\gamma} g_{i,\gamma} y_{\alpha+\beta+\gamma})_{|\alpha|,|\beta|\leqslant d-k_i}\\[.2cm]
k_i := \max\{ \lceil|\alpha|/ 2 \rceil ~\text{s.t.}~ g_{i,\alpha} \neq 0 \}.
\end{array}
\end{equation}
Above, $\lceil \cdot \rceil$ denotes the ceiling of a real number. Note that the relaxation order~$d$ must be greater than or equal to $d^\text{min} := \max_i k_i$. See Section~\ref{sec:Numerical experiments} for an illustrative example of a second-order moment/sum-of-squares relaxation of a quadratically constrained quadratic program corresponding to a two-bus OPF problem. 

%{\color{red}Consider, for instance, the following polynomial optimization problem corresponding to the second-order moment/sum-of-squares relaxation for a two-bus OPF problem:}

Our approach extends a relaxation technique proposed in~\cite{hijazi-pscc2016} that leverages theory developed in~\cite{lasserre2010representations,lasserre2011convex} regarding the importance of the convexity of the feasible space rather than its algebraic representation.\footnote{In particular, under a mild degeneracy condition, log-barrier solution methods are guaranteed to converge to a global minimizer of a problem with a convex feasible space even if the algebraic representation of the feasible space is non-convex.} The technique in~\cite{hijazi-pscc2016} further relaxes the Shor relaxation~\cite{shor-1987b} using necessary conditions for  positive semidefiniteness based on matrix determinants. Specifically, the approach in~\cite{hijazi-pscc2016} exploits the fact that a symmetric matrix is positive semidefinite if and only if all its principal minors (i.e., the determinants of its principal submatrices) are non-negative~\cite{prussing1986principal}. Extending this approach to the Lasserre hierarchy~\eqref{eq:lasserre} yields
\begin{equation}
\label{eq:lasserre_det}
\begin{array}{ll}
\inf_y & L_y(f)  \\
\mathrm{s.t.} & y_0 = 1 \\
 & \det_I(M_d(y)) \geqslant 0 \\
 & \quad \text{for all~} I = \{1\}, \{2\}, \ldots, \{1,2\}, \\ & \{1,3\}, \ldots,  \{1,2,\ldots,n_0\}, \\
 & \det_{J_i}(M_{d-k_i}(g_iy)) \geqslant 0 \\
 & \quad \text{for all~} J_i = \{1\}, \{2\}, \ldots, \{1,2\},\\ & \{1,3\}, \ldots, \{1,2,\ldots,n_i\}, \\
 & \qquad i=1,\ldots,m,
\end{array}
\end{equation}
where $n_0 := \frac{\left(n+d\right)!}{n!\,d!}$ is the size of the matrix $M_d(y)$, $n_i := \frac{\left(n+d-k_i\right)!}{n!\,\left(d-k_i\right)!}$ is the size of the matrix $M_{d-k_i}(y)$, and $\det_I$ denotes the principal minor corresponding to the index $I$ (i.e., the determinant of the submatrix whose rows and columns are indexed by $I$), likewise for $\det_{J_i}$.

To derive our approach, we exploit the KKT conditions of the formulation~\eqref{eq:lasserre_det}. To express these conditions, let $\lambda$ denote the KKT multiplier associated with the constraint $y_0 = 1$. Let $\lambda_{0,I}$ denote the KKT multiplier associated with the constraint $\det_I(M_d(y)) \geqslant 0$. Let $\lambda_{i,J_i}$ denote the KKT multiplier associated with the constraint $\det_{J_i}(M_{d-k_i}(g_iy)) \geqslant 0$. The KKT conditions of~\eqref{eq:lasserre_det} are
\begin{equation}
\begin{array}{l}
\nabla L_y(f) - \lambda e_1 - \sum\limits_I \lambda_{0,I} \nabla_y \det_I(M_d(y)) + \hdots \\ - \sum\limits_{i=1}^m \sum\limits_{J_i} \lambda_{i,J_i} \nabla_y \det_{J_i}(M_{d-k_i}(g_iy)) = 0, \\
y_0 = 1,\\
\det_I(M_d(y)) \geqslant 0, \\
\det_{J_i}(M_{d-k_i}(g_iy)) \geqslant 0, \quad i=1,\ldots,m, \\
\lambda_{0,I} \geqslant 0, \\
\lambda_{i,J_i} \geqslant 0, \quad i=1,\ldots,m, \\
\det_I(M_d(y)) \lambda_{0,I} = 0, \\
\det_{J_i}(M_{d-k_i}(g_iy)) \lambda_{i,J_i} = 0, \quad i=1,\ldots,m.
\end{array}
\end{equation}

We have that $\nabla \det (M) = \text{co}(M) $ where the comatrix $\text{co}(M)$ is defined by $\text{co}(M)_{i,j} := (-1)^{i+j}\det(M^{(i,j)})$; $M^{(i,j)}$ is the matrix $M$ where row $i$ and column $j$ have been removed. Let $\text{co}_I(\cdot)$ denote the comatrix of the submatrix whose rows and columns are indexed by a subset of indices $I$. Below, we will consider matrices $B_\alpha$ and $C_{i,\alpha}$ such that the moment and localizing matrices are as follows: 
\begin{equation}
M_d(y) =: \sum_{|\alpha| \leqslant 2d} B_{\alpha} y_{\alpha}
\end{equation}
and
\begin{equation}
M_{d-k_i}(g_iy) =: \sum_{|\alpha| \leqslant 2(d-k_i)} C_{i,\alpha}\, y_{\alpha}.
\end{equation}
 The KKT conditions are then\footnote{We use the fact that for all sequences of multi-indexed real numbers $(h_\alpha)_{|\alpha|\leqslant 2d}$, we have $d(\det \circ M_d)(y).h = d(\det)(M_d(y))(d(M_d)(y).h) = d(\det)(M_d(y)).M_d(h) = \text{trace}[\text{co}(M_d(y)) M_d(h)] = \text{trace}[\text{co}(M_d(y)) \sum_\alpha B_\alpha h_\alpha] = \sum\limits_\alpha \text{trace}[\text{co}(M_d(y)) B_\alpha] h_\alpha$.}

\begin{equation}
\label{eq:nonlin_feasibility}
\begin{array}{l}
f_{\alpha} - \lambda e_1 - \sum\limits_I \lambda_{0,I}\, \text{trace}(\text{co}_I(M_d(y)) B_{I,\alpha})+\hdots \\ - \sum\limits_{i=1}^m \sum\limits_{J_i} \lambda_{i,J_i} \text{trace}(\text{co}_{J_i}(M_{d-k_i}(g_iy)) C_{i,J_i,\alpha}) = 0, \\
y_0 = 1,\\
\det_I(M_d(y)) \geqslant 0, \\
\det_{J_i}(M_{d-k_i}(g_iy)) \geqslant 0, \quad i=1,\ldots,m, \\
\lambda_{0,I} \geqslant 0, \\
\lambda_{i,J_i} \geqslant 0, \quad i=1,\ldots,m, \\
\det_I(M_d(y)) \lambda_{0,I} = 0, \\
\det_{J_i}(M_{d-k_i}(g_iy)) \lambda_{i,J_i} = 0, \quad i=1,\ldots,m.
\end{array}
\end{equation}

If we substitute $y$ using the feasible point $\hat{x}$, i.e., by taking $\hat{y} := ( \hat{x}^\alpha )_{|\alpha|\leqslant 2d}$, the KKT conditions reduce to a linear feasibility problem, that is
\begin{equation}
\label{eq:lp_feasibility}
\begin{array}{l}
f_{\alpha} - \lambda e_1 - \sum\limits_I \lambda_{0,I} \text{trace}(\text{co}_I(M_d(\hat{y})) B_{I,\alpha}) + \hdots \\ - \sum\limits_{i=1}^m \sum\limits_{J_i} \lambda_{i,J_i} \text{trace}(\text{co}_{J_i}(M_{d-k_i}(g_i\hat{y})) C_{i,J_i,\alpha}) = 0, \\
    \lambda_{0,I} \geqslant 0, \\
\lambda_{i,J_i} \geqslant 0, \quad i=1,\ldots,m, \\
\hat{y}_{2\alpha}\, \lambda_{0,\{\alpha\}} = 0, \quad \forall |\alpha| \leqslant d, \\[.1cm]
(\sum\limits_{\beta} g_{i,\beta}\, \hat{y}_{2\alpha+\beta})\, \lambda_{i,\{\alpha\}} = 0, \quad \forall |\alpha| \leqslant d-k_i, \\
i=1,\ldots,m.
\end{array}
\end{equation}

Observe that all the determinants are equal to zero, except when the index sets $I$ and $J_i$ are singletons (i.e.,  $I=J_i=\{\alpha\}$ where $\alpha \in \mathbb{N}^n$). In that case, the determinant is simply equal to a diagonal term of the matrix. 

Any point that satisfies~\eqref{eq:lp_feasibility} provides values for dual variables that, in combination with the primal variables $\hat{y}$ obtained from the candidate solution $\hat{x}$, satisfy the KKT conditions corresponding to a convex relaxation of~\eqref{eq:pop}. Hence, any feasible point for~\eqref{eq:lp_feasibility} certifies global optimality of the candidate feasible point $\hat{x}$ for the polynomial optimization problem~\eqref{eq:pop}.\footnote{Observe that neither \cite[Assumption 2.1]{lasserre2010representations} nor Slater's condition are required in order to show that every KKT point is a minimizer in \cite[Theorem 2.3]{lasserre2010representations}. In other words, all assumptions needed for the theory in \cite{lasserre2010representations} to be applicable for our purposes are satisfied in our setting.}

Since all values $\hat{y}$ are fixed, we emphasize that~\eqref{eq:lp_feasibility} is \emph{linear} in the KKT multipliers $\lambda$, $\lambda_{0,I}$, and $\lambda_{0,J_i}$. Thus, the KKT conditions~\eqref{eq:lp_feasibility} can be verified using a linear programming solver, which is significantly more computationally tractable than solving the KKT conditions for the relaxation~\eqref{eq:nonlin_feasibility} with both primal variables $y$ and dual variables $\lambda$, $\lambda_{0,I}$, and $\lambda_{0,J_i}$ allowed to vary (a complementarity problem).

\section{Numerical experiments}
\label{sec:Numerical experiments}

%{\color{red}Exploit sparsity and selectively apply the higher-order constraints. Could also use complex hierarchy rather than real hierarchy. Most importantly, choose where to enforce the higher-order constraints based on where the slack variables are non-zero.}

We implemented the algorithms with Matlab R2020a. We used Gloptipoly 3 to obtain the Lasserre hierarchy formulation. As for solving the KKT conditions, we used the LINPROG and LSQLIN solver provided by Matlab. The optimality tolerance of both solvers are set to $10^{-8}$.

\subsection{Univariate example}

\begin{figure}
    \centering
    \includegraphics[width=.45\textwidth]{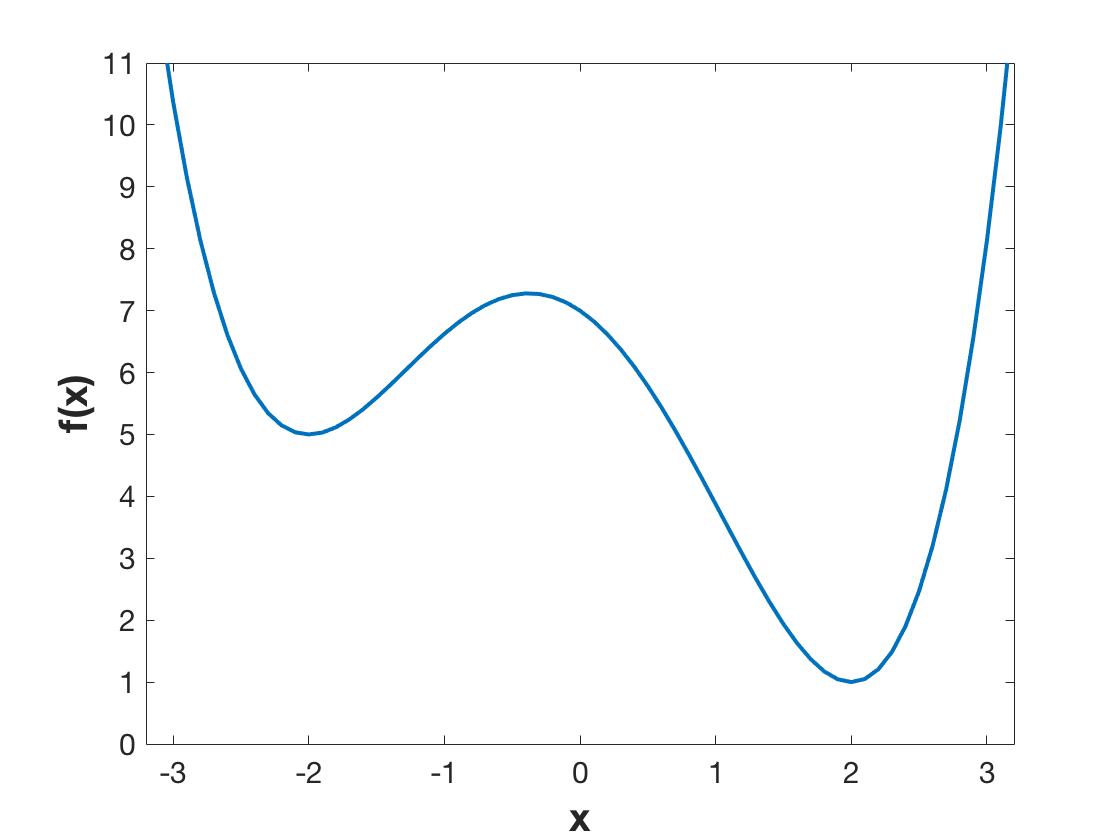}
    \caption{Univariate example}
    \label{fig:uni}
\end{figure}

Consider the following polynomial optimization problem represented in Figure \ref{fig:uni} with a local minimum at $-2$ and a global minimum at $2$:
\begin{equation}
    \inf_{x\in \mathbb{R}} ~~~ \frac{1}{4}x^4+\frac{1}{8}x^3-2x^2-\frac{3}{2}x+7 ~~~\text{s.t.} ~~~ - x^2 + 5 \geqslant 0.
\end{equation}

% $K=\{x\in\mathbb{R}\,|\, 5 - x^2\geq 0\}$

% Riesz function: $L_y(f) = \frac{1}{4}y_4+\frac{1}{8}y_3-2y_2-\frac{3}{2}y_1+7y_0$

% Moment matrix: 
% $M_2(y) = (y_{\alpha+\beta})_{0\leq \alpha, \beta\leq 2}\succeq0$

% Localizing matrix: $M_1(gy)=(-y_{\alpha+\beta+2}+5y_{\alpha+\beta})_{0\leq \alpha, \beta\leq 1} \succeq 0$

The second-order moment relaxation reads:
\begin{equation}
\begin{array}{ll}
     \inf_{y \in \mathbb{R}^5} & \frac{1}{4}y_4+\frac{1}{8}y_3-2y_2-\frac{3}{2}y_1+7y_0 \\[.4cm]
     \text{s.t.} &
     \left\{
     \begin{array}{l}
          y_0 = 1, \\
          \begin{pmatrix}
            y_0 & y_1 & y_2 \\
            y_1 & y_2 & y_3 \\
            y_2 & y_3 & y_4 \\
          \end{pmatrix} \succcurlyeq 0, \\[1.5em]
          \begin{pmatrix}
            -y_2 + 5y_0 & -y_3 + 5y_1 \\
            -y_3 + 5y_1 & -y_4 + 5y_2 \\
          \end{pmatrix} \succcurlyeq 0.
     \end{array}
     \right.
\end{array}
\end{equation}

Using the determinant condition for positive semidefiniteness as in~\eqref{eq:lasserre_det}, an equivalent formulation of the second-order moment relaxation is as follows:
\begin{equation}
\label{eq:non-convex1}
     \inf_{y \in \mathbb{R}^5} ~~~ \frac{1}{4}y_4+\frac{1}{8}y_3-2y_2-\frac{3}{2}y_1+7y_0
\end{equation}
subject to
\begin{equation}
\label{eq:non-convex2}
     \left\{
     \begin{array}{l}
          y_0 = 1,      \\[1mm]                                                          															 
          y_0 \geqslant 0,   ~                                                													  
          y_2 \geqslant 0,    ~                                                 														
          y_4 \geqslant 0, 						\\[1mm]																				
         y_0y_2 - y_1^2 \geqslant 0, ~ 																	
         y_0y_4 - y_2^2 \geqslant 0, ~ 																				
         y_2y_4 - y_3^2  \geqslant 0, 									\\[1mm]													 
         y_0(y_2y_4-y_3^2) - y_1(y_1y_4-y_3y_2) \\ + y_2(y_1y_3-y_2^2) \geqslant 0,		\\[1mm]
         -y_2 + 5 y_0 \geqslant 0, 	~																												 
         -y_4 + 5y_2 \geqslant 0, 										\\[1mm]																		
         (-y_2+5y_0)(-y_4+5y_2) - (-y_3+5y_1)^2 \geqslant 0. 																	 
     \end{array}
     \right.
\end{equation}
Observe that each constraint taken separately may be non-convex, yet together they yield a convex region. We then define the following Lagrangian:
\begin{equation}
\begin{array}{cc}
     \mathcal{L}(y,\lambda,\lambda_0,\lambda_1) := \frac{1}{4}y_4+\frac{1}{8}y_3-2y_2-\frac{3}{2}y_1 \\ +7y_0 + (1-y_0){\lambda}  - y_0 {\lambda_{0,\{1\}}} - y_2 {\lambda_{0,\{2\}}} \\- y_4{\lambda_{0,\{3\}}} 
     -(y_0y_2 - y_1^2){\lambda_{0,\{1,2\}}} 
         \\-(y_0y_4 - y_2^2) {\lambda_{0,\{1,3\}}} -  
         (y_2y_4 - y_3^2){\lambda_{0,\{2,3\}}} \\   
         - [
         y_0(y_2y_4-y_3^2) - y_1(y_1y_4-y_3y_2) \\+ y_2(y_1y_3-y_2^2)]{\lambda_{0,\{1,2,3\}}} +  \\
         -(-y_2 + 5 y_0) {\lambda_{1,\{1\}}} -(-y_4 + 5y_2){\lambda_{1,\{2\}}} + \\
         -[(-y_2+5y_0)(-y_4+5y_2) \\- (-y_3+5y_1)^2]{\lambda_{1,\{1,2\}}}.
\end{array}
\end{equation}

The KKT conditions that need to be met by the dual multipliers of \eqref{eq:non-convex1}--\eqref{eq:non-convex2} are thus given by:

\begin{equation}
\left\{    
    \begin{array}{l}
       7 - {\lambda} - {\lambda_{0,\{1\}}} -y_2{\lambda_{0,\{1,2\}}} -y_4{\lambda_{0,\{1,3\}}} \\ -(y_2y_4-y_3^2){\lambda_{0,\{1,2,3\}}} = 0, \\[.2cm]
       -3/2 + 2y_1{\lambda_{0,\{1,2\}}} \\- (-2y_1y_4 + 2y_3y_2){\lambda_{0,\{1,2,3\}}}\\ + 10(-y_3+5y_1){\lambda_{1,\{1,2\}}}= 0, \\[.2cm]
         -2 - {\lambda_{0,\{2\}}} 
     -y_0{\lambda_{0,\{1,2\}}} +
         2y_2^2 {\lambda_{0,\{1,3\}}} \\- 
         y_4{\lambda_{0,\{2,3\}}}  
         +(- y_0y_4 - y_1y_3 - y_1y_3 \\ + 3y_2^2) 
         {\lambda_{0,\{1,2,3\}}} +
          {\lambda_{1,\{1\}}} - 5{\lambda_{1,\{2\}}} \\
          +(-y_4 + 10y_2 -25y_0){\lambda_{1,\{1,2\}}} = 0, \\[.4cm]
       1/8 + 2y_3{\lambda_{0,\{2,3\}}} \\- (-2y_0y_3+ 2y_1y_2) {\lambda_{0,\{1,2,3\}}} \\ - (-2y_3 + 10y_1){\lambda_{1,\{1,2\}}} = 0, \\[.4cm]
       1/4 - {\lambda_{0,\{3\}}}
         - y_0 {\lambda_{0,\{1,3\}}} - 
         y_2{\lambda_{0,\{2,3\}}} \\
         - (y_0y_2 - y_1^2) {\lambda_{0,\{1,2,3\}}} +{\lambda_{1,\{2\}}}\\ + (-y_2+5y_0){\lambda_{1,\{1,2\}}}= 0, \\[.2cm]
        %  y_0 = 1,      \\[1mm]                                                          															 
        %   y_0 \geqslant 0,   ~                                                													  
        %   y_2 \geqslant 0,    ~                                                 														
        %   y_4 \geqslant 0, 						\\[1mm]																				
        %  y_0y_2 - y_1^2 \geqslant 0, ~ 																	
        %  y_0y_4 - y_2^2 \geqslant 0, ~ 																				
        %  y_2y_4 - y_3^2  \geqslant 0, 									\\[1mm]													 
        %  y_0(y_2y_4-y_3^2) - y_1(y_1y_4-y_3y_2) + y_2(y_1y_3-y_2^2) \geqslant 0,		\\[1mm]
        %  -y_2 + 5 y_0 \geqslant 0, 	~																												 
        %  -y_4 + 5y_2 \geqslant 0, 										\\[1mm]																		
        %  (-y_2+5y_0)(-y_4+5y_2) - (-y_3+5y_1)^2 \geqslant 0, \\[2mm] 		
         {\lambda_{0,\{1\}}},~{\lambda_{0,\{2\}}},~{\lambda_{0,\{3\}}},~ {\lambda_{0,\{1,2\}}},\\{\lambda_{0,\{1,3\}}},~{\lambda_{0,\{2,3\}}},~{\lambda_{0,\{1,2,3\}}} \geqslant 0, \\
         {\lambda_{1,\{1\}}},~{\lambda_{1,\{2\}}},~ {\lambda_{1,\{1,2\}}} \geqslant 0, \\[.2cm]
         y_0 {\lambda_{0,\{1\}}} = y_2 {\lambda_{0,\{2\}}} = y_4 {\lambda_{0,\{3\}}} = 0, \\[.2cm]
         {\lambda_{1,\{1\}}}(-y_2+5y_0) = {\lambda_{1,\{2\}}}(-y_4+5y_2) = 0.
         %\textcolor{red}{\text{Be sure to add other KKT conditions, i.e primal-dual feasibility and complementary slackness}}
    \end{array}
\right.
\end{equation}

We next check whether $x=-2$ or $x=2$ are global solutions. In order to do so, we plug $x$ into the primal variable $y$ in the above system, i.e., $y_i := x^i$, and see whether there exists a feasible set of dual variables. If there exists feasible dual variables, then the point $x$ is a global solution. If not, then we may need to increase the order of the hierarchy and try again. In order to check for linear feasibility, we respectively minimize the $\ell_1$ and $\ell_2$ norms of the equality constraints, subject to the non-negativity constraints of the dual multipliers. The $\ell_1$ norm results in a linear program while the $\ell_2$ norm results in a convex quadratic optimization problem. The numerical results are shown in Table~\ref{tab:uni_example}. The local solution yields an infeasible linear system, while the global solution yields a linear system that is feasible with very high accuracy.

% To solve the KKT system, we solve the following problem:
% \begin{equation}
%     \begin{split}
%         \min\; & \|C \lambda - d\| \\
%         \text{s.t.}\; & \lambda_{-1} \geq 0 \\
%     \end{split}
% \end{equation}
% where $\lambda_{-1}$ means all the $\lambda$ except for the one for the $y_0=1$ constraint. 
%We tested the KKT system on two solutions of the univariate example, and 
% \begin{table}[h!]
% \centering
% \label{tab:uni_example}
%     \begin{tabular}{|c|c|c|c|}
%     \hline 
%     Solution & Optimality & Residual in $\ell_1$-norm & Residual in $\ell_2$-norm \\
%     \hline 
%     $x=-2$ & Local & $1.55e+00$ & $2.25e+00$ \\
%     $x=\hphantom{-}2$ & Global & $2.30e-07$ & $1.48e-08$ \\
%     \hline
% \end{tabular}
% \caption{Global optimality check in univariate example}
% \end{table}

\begin{table}[h!]
\caption{Global optimality check in univariate example}
\label{tab:uni_example}
\begin{center}
\begin{tabular}{|c|c|c|c|}
    \hline 
    Solution & Optimality & $\ell_1$ norm & $\ell_2$ norm \\
    \hline 
    $x=-2$ & Local & $1.55\times 10^{0\hphantom{-}}$ & $2.25\times 10^{0\hphantom{-}}$ \\
    $x=\hphantom{-}2$ & Global & $2.30\times 10^{-7}$ & $1.48\times 10^{-8}$ \\
    \hline
\end{tabular}
\end{center}
\end{table}

\subsection{Bivariate example}

\begin{figure}
    \centering
    \includegraphics[width=.45\textwidth]{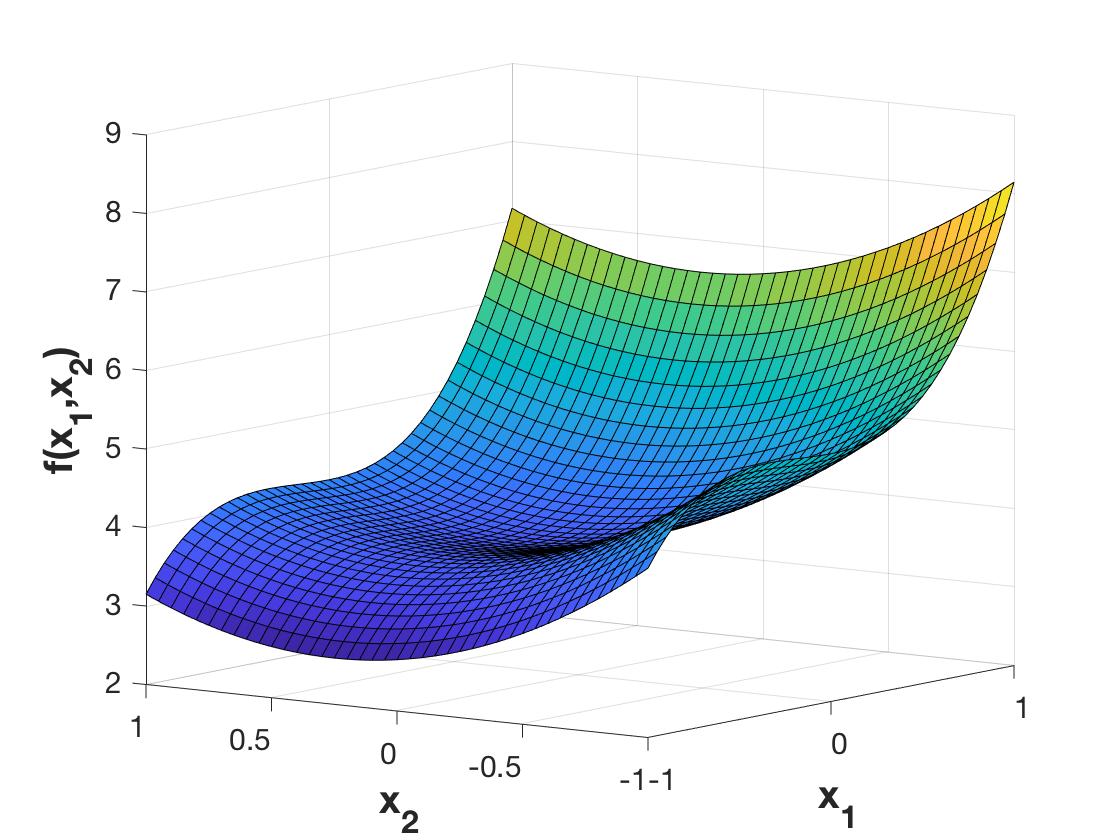}
    \caption{Bivariate example}
    \label{fig:bivariate}
\end{figure}

Consider the following polynomial optimization problem represented in Figure \ref{fig:bivariate}:
\begin{align}\nonumber
        & \inf_{x_1,x_2\in \mathbb{R}} ~~~ 2x_1^3 + x_1^2 + \frac{1}{4}x_1x_2 + x_2^2 -\frac{1}{2} x_2 +\frac{1}{16} \\
        & \text{subject to} ~~ x_1^2 + x_2^2 \leq 1.
\end{align}

The numerical results in Table~\ref{tab:biv_example} again show that the method correctly distinguishes the global solution from the local solution. The KKT system stemming from the second-order moment relaxation was used.
% \begin{center}
% \label{tab:biv_example}
%     \begin{tabular}{|c|c|c|c|}
%     \hline 
%     Solution & Optimality & L1 Residual Norm & L2 Residual Norm \\
%     \hline 
%     $x_1=-0.0360, ~x_2=0.254$ & Local & $1.991996e+00$ & $3.684069e+00$ \\
%     $x_1=-0.9922, ~x_2=0.125$ & Global & $1.747277e-04$ & $7.133155e-08$ \\
%     \hline
% \end{tabular}
% \end{center}
% \begin{table}[h!]
% \centering
% \label{tab:biv_example}
%     \begin{tabular}{|c|c|c|c|}
%     \hline 
%     Solution & Optimality & Residual in $\ell_1$-norm & Residual in $\ell_2$-norm \\
%     \hline 
%     $x_1=-0.036, ~x_2=0.254$ & Local & $1.99e+00$ & $3.68e+00$ \\
%     $x_1=-0.992, ~x_2=0.125$ & Global & $1.75e-04$ & $7.13e-08$ \\
%     \hline
% \end{tabular}
% \caption{Global optimality check in bivariate example}
% \end{table}

\begin{table}[h!]
\caption{Global optimality check in bivariate example}
\vspace*{-1.5em}
\label{tab:biv_example}
\begin{center}
\begin{tabular}{|c|c|c|c|}
    \hline 
    Solution & Optimality & $\ell_1$ norm & $\ell_2$ norm \\
    \hline 
    $x=(-0.036,0.254)$ & Local & $1.99\times 10^{0\hphantom{-}}$ & $3.68\times 10^{0\hphantom{-}}$ \\
    $x=(-0.992,0.125)$ & Global & $1.75\times 10^{-4}$ & $7.13\times 10^{-8}$ \\
    \hline
\end{tabular}
\end{center}
\end{table}

\subsection{Trivariate example}

\begin{figure}[h!]
  \centering
    \includegraphics[width=.49\textwidth]{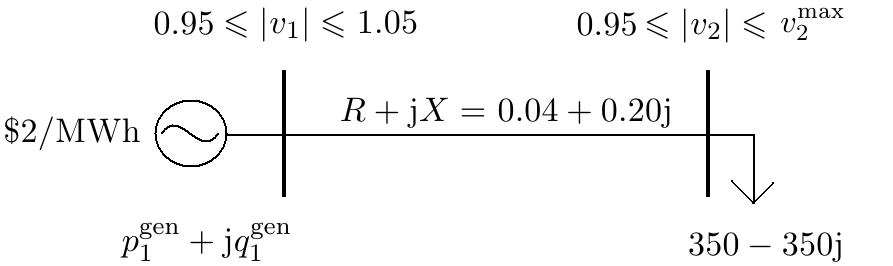}
  \caption{WB2 2-bus system}
  \label{fig:WB2 2-bus system}
\end{figure}

Consider the following polynomial optimization problem:
\begin{equation}
        \inf_{x_1,x_2,x_3 \in \mathbb{R}} \frac{2500}{13}x_1^2 - \frac{12500}{13}x_1x_3 - \frac{2500}{13}x_1x_2 
        \end{equation}
        subject to
        \begin{equation}
        \left\{
        \begin{array}{c}
              \frac{25}{26}x_1x_2 - \frac{125}{26}x_1x_3 - \frac{25}{26}x_2^2 - \frac{25}{26}x_3^2 - \frac{7}{2} = 0 \\[2mm]
             \frac{125}{26}x_1x_2 + \frac{25}{26}x_1x_3 - \frac{125}{26}x_2^2 - \frac{125}{26}x_3^2 + \frac{7}{2} = 0 \\[2mm]
             0 \leqslant \frac{25}{26}x_1^2 - \frac{125}{26}x_1x_3 - \frac{25}{26}x_1x_2 \leqslant 6 \\[2mm]
             -4 \leqslant \frac{25}{26}x_1x_3 - \frac{125}{26}x_1x_2 + \frac{125}{26}x_1^2 \leqslant 4 \\[2mm]
             0.9025 \leqslant x_1^2 \leqslant 1.1025 \\[2mm]
        0.9025 \leqslant x_2^2 + x_3^2 \leqslant 1.1025
        \end{array}
        \right.
\end{equation}
% variables: Vd1, Vd2, Vq2;
% objective: 
% (2500.*Vd1.^2)./13 - (12500.*Vd1.*Vq2)./13 - (2500.*Vd1.*Vd2)./13;
% constraints: 
% (25.*Vd1.*Vd2)./26 - (125.*Vd1.*Vq2)./26 - (25.*Vd2.^2)./26 - (25.*Vq2.^2)./26 - 7./2 = 0;
% (125.*Vd1.*Vd2)./26 + (25.*Vd1.*Vq2)./26 - (125.*Vd2.^2)./26 - (125.*Vq2.^2)./26 + 7./2 = 0;
% 0 <= (25.*Vd1.^2)./26 - (125.*Vd1.*Vq2)./26 - (25.*Vd1.*Vd2)./26 <= 6;
% -4 <= (25.*Vd1.*Vq2)./26 - (125.*Vd1.*Vd2)./26 + (125.*Vd1.^2)./26 <= 4;
% 0.9025 <= Vd1.^2 <= 1.1025;
% 0.9025 <= Vd2.^2 + Vq2.^2 <= 1.1025;
% solution 1: 
% Vd1 = 0.95, Vd2 = 0.41343, Vq2 = -0.88419;
% solution 2: 
% Vd1 = 0.95234, Vd2 = 0.56965, Vq2 = -0.88204;
This problem corresponds to the two-bus alternating current optimal power flow problem in Figure~\ref{fig:WB2 2-bus system} from~\cite{bukhsh2013}.

\begin{figure}[h!]
    \centering
    \includegraphics[width=.49\textwidth]{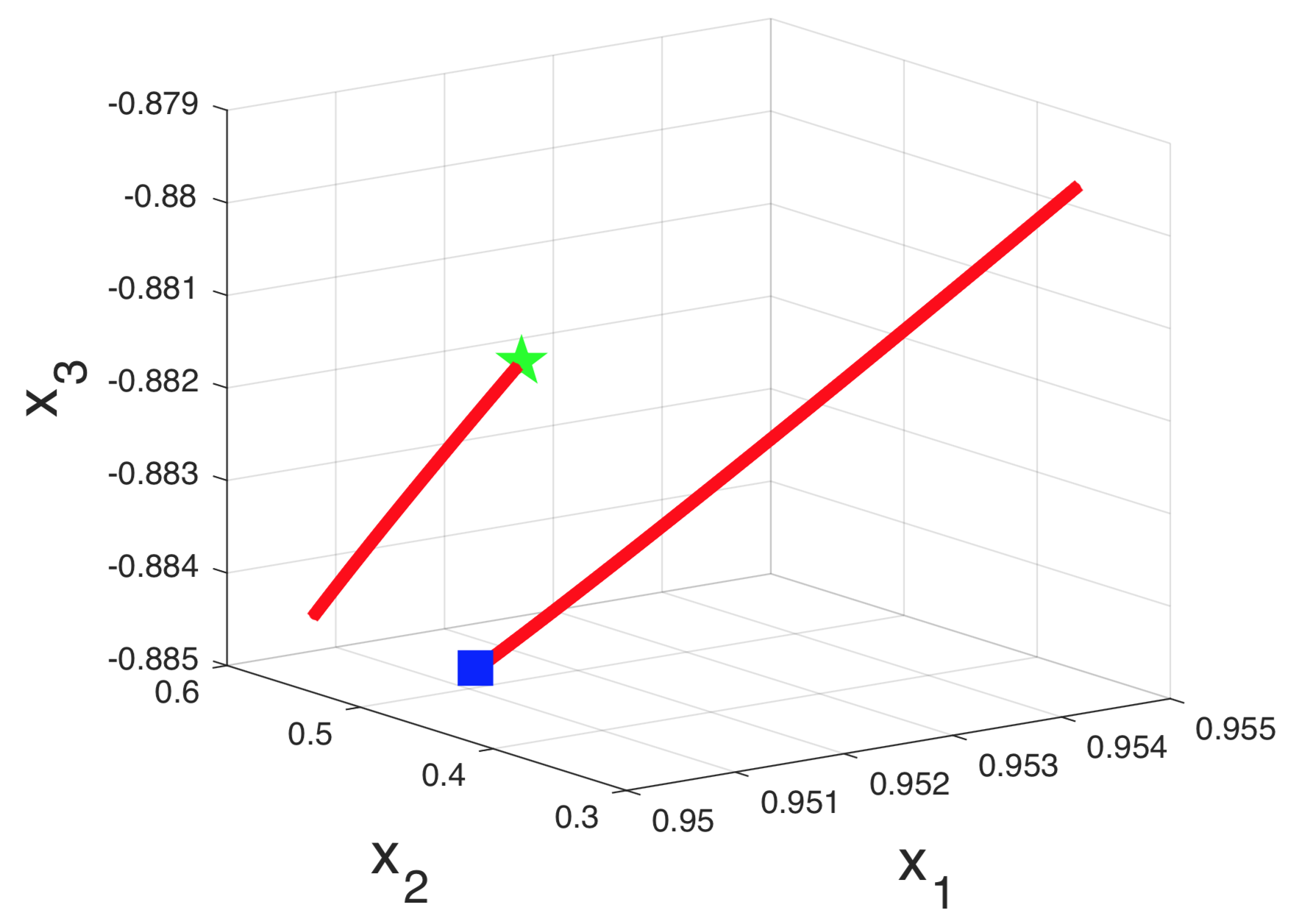}
    \caption{Trivariate example}
    \label{fig:trivariate}
\end{figure}

The feasible region of the polynomial optimization problem is represented in Figure \ref{fig:trivariate}, along with a local solution (blue square) and a global solution (green star).

The numerical results in Table \ref{tab:tri_example} show that the method is able to distinguish between the local and global solutions, albeit less strikingly than in the univariate and bivariate examples.  We attribute this to the proximity between the local and global solutions, which are taken from the test case archive.\footnote{\href{https://www.maths.ed.ac.uk/optenergy/LocalOpt/WB2.html}{\url{https://www.maths.ed.ac.uk/optenergy/LocalOpt/WB2.html}}} Their objective values are also close: the objective values of the local and global solutions are $905.73$ and $877.78$, respectively (a 3.2\% difference). We used the KKT system stemming from the second-order moment relaxation, consistent with the empirical finding in~\cite{cedric_tps,pscc2014} that this relaxation is often tight in the context of power systems.

% \begin{table}[h!]
% \label{tab:tri_example}
% \centering
%     \begin{tabular}{|c|c|c|c|}
%     \hline 
%     Solution & Optimality & Res. in $\ell_1$-norm & Res. in $\ell_2$-norm \\
%     \hline 
%     $x_1=0.95000, x_2=0.41343, x_3=-0.88419$ & Local & $7.21e-02$ & $2.96e-03$ \\
%     %$x_1=-0.95000,x_2= -0.41343, x_3=0.88419$ & Local & $7.206157e-02$ & $2.961877e-03$ \\
%     $x_1=0.95234, x_2=0.56965, x_3=-0.88204$ & Global & $5.31e-03$ & $9.93e-06$ \\
%     %$x_1=-0.95234, x_2=-0.56965, x_3=0.88204$ & Global & $5.307450e-03$ & $9.930022e-06$ \\
%     \hline
% \end{tabular}
% \caption{Global optimality check in trivariate example}
% \end{table}

\begin{table}[h!]
\caption{Global optimality check in trivariate example}
\label{tab:tri_example}
\vspace*{-2em}
\begin{center}
\begin{tabular}{|c|c|c|c|}
    \hline 
    Solution & Optimality & $\ell_1$ norm & $\ell_2$ norm \\
    \hline 
    %$x=(0.95000,0.41343,-0.88419)$ 
    $x=(.950,.413,-.884)$ 
    & Local & $7.21\times 10^{-2}$ & $2.96\times 10^{-3}$ \\
    %$x_1=-0.95000,x_2= -0.41343, x_3=0.88419$ & Local & $7.206157e-02$ & $2.961877e-03$ \\
    %$x=(0.95234,0.56965,-0.88204)$ 
    $x=(.952,.570,-.882)$
    & Global & $5.31\times 10^{-3}$ & $9.93\times 10^{-6}$ \\
    %$x_1=-0.95234, x_2=-0.56965, x_3=0.88204$ & Global & $5.307450e-03$ & $9.930022e-06$ \\
    \hline
\end{tabular}
\end{center}
\end{table}

We conclude this section with the runtimes in Table \ref{tab:my_label}. The runtimes are the average of 5 repeated experiments. One can see that there is a significant speed-up, by at least an order of magnitude, when checking for global optimality instead of computing a global solution. We do not report the time taken to setup the linear program and the semidefinite program, both of which could be reduced with a more sophisticated implementation.
\begin{table}[]
    \caption{Runtimes (in seconds)}
    \centering
    \begin{tabular}{|c|c|c|c|}
    \hline 
        Problem & SDP & $\ell_1$ norm & $\ell_2$ norm \\
    \hline 
        Univariate & $0.960$ & $0.042$ & $0.044$ \\
        Bivariate & $0.922$ & $0.032$ & $0.042$ \\
        Trivariate & $2.068$ & $0.050$ & $0.706$ \\
    \hline 
    \end{tabular}
    \label{tab:my_label}
\end{table}

\section{Conclusion}
\label{sec:conclusion}

We have proposed a method for verifying global optimality by checking for the feasibility of a linear system of equations. We have demonstrated the proof of concept on several small examples of polynomial optimization problems.
A topic of further investigation is to exploit sparsity in order to make the approach tractable on larger optimal power flow instances. One could use techniques involving chordal graphs as was done in \cite{josz2018lasserre}, which enabled the application of the Lasserre hierarchy to large-scale instances. Since global solutions were extracted in many instances within minutes, it is reasonable to believe that checking for global optimality should be possible with significantly faster runtimes. Another possible avenue for future research is to find another equivalent formulation of the Lasserre hierarchy whose KKT conditions would also lead to a tractable system for the dual multipliers (perhaps a second-order conic program).

\bibliography{mybib}{}
\bibliographystyle{ieeetr}

% \begin{thebibliography}{99}

% \bibitem{c1}
% J.G.F. Francis, The QR Transformation I, {\it Comput. J.}, vol. 4, 1961, pp 265-271.

% \bibitem{c2}
% H. Kwakernaak and R. Sivan, {\it Modern Signals and Systems}, Prentice Hall, Englewood Cliffs, NJ; 1991.

% \bibitem{c3}
% D. Boley and R. Maier, "A Parallel QR Algorithm for the Non-Symmetric Eigenvalue Algorithm", {\it in Third SIAM Conference on Applied Linear Algebra}, Madison, WI, 1988, pp. A20.

% \end{thebibliography}

\end{document}